\documentclass[aop,seceqn,nameyear,MSNbibl,dvips]{arximspdf}

\doi{10.1214/10-AOP530}
\volume{38}
\issue{5}
\pubyear{2010}
\firstpage{2009}
\lastpage{2022}

\makeatletter

\newtheorem{lemma}{Lemma}[section]
\newtheorem{proposition}{Proposition}[section]
\newtheorem{theorem}{Theorem}[section]
\newproclaim{remark}{Remark}[section]
\newproclaim{example}{Example}[section]

\makeatother

\begin{document}
\begin{frontmatter}

\title{Central limit theorem for Fourier transforms of
stationary processes}
\runtitle{CLT for Fourier transforms}

\begin{aug}
\author[A]{\fnms{Magda} \snm{Peligrad}\thanksref{T1}\ead[label=e1]{peligrm@math.uc.edu}} and
\author[B]{\fnms{Wei Biao} \snm{Wu}\corref{}\thanksref{T2}\ead[label=e2]{wbwu@galton.uchicago.edu}}
\runauthor{M. Peligrad and W. B. Wu}
\affiliation{University of Cincinnati and University of Chicago}
\address[A]{Department of Mathematical Sciences\\
University of Cincinnati\\
PO Box 210025\\
Cincinnati, Ohio 45221\\
USA\\
\printead{e1}}
\address[B]{Department of Statistics\\
University of Chicago\\
5734 S. University Avenue\\
Chicago, Illinois 60637\\
USA\\
\printead{e2}}
\end{aug}

\thankstext{T1}{Supported in part by a Charles Phelps
Taft Memorial Fund grant and NSA Grant H98230-09-1-0005.}

\thankstext{T2}{Supported in part by a Charles Phelps Taft
Memorial Fund and NSF Grants DMS-04-78704 and
0906073.}

\received{\smonth{6} \syear{2009}}
\revised{\smonth{12} \syear{2009}}

%
\begin{abstract}
We consider asymptotic behavior of Fourier transforms of
stationary ergodic sequences with finite second moments. We
establish a central limit theorem (CLT) for almost all frequencies
and also an annealed CLT. The theorems hold for all regular
sequences. Our results shed new light on the foundation of
spectral analysis and on the asymptotic distribution of
periodogram, and it provides a nice blend of harmonic analysis,
theory of stationary processes and theory of martingales.
\end{abstract}

%

\end{frontmatter}

\section{Introduction}\label{sec1}
In frequency or spectral domain analysis of time series,
periodograms play a fundamental role. Since its introduction by
\citet{Schuster1898}, periodograms have been used in almost all
scientific fields. Given a realization $(X_{j})_{j=1}^{n}$ of a
stochastic process $(X_{j})_{j\in\mathbb{Z}}$, the periodogram is
defined as
\[
I_{n}(\theta)={\frac{1}{{2\pi n}}} \Biggl\vert
\sum_{j=1}^{n}X_{j}\exp(ji\theta) \Biggr\vert
^{2}, \qquad\theta\in\mathbb{R},
\]
where $i=\sqrt{-1}$ is the imaginary unit. Periodogram is the
building block in spectral domain analysis and a distributional
theory is clearly needed in the related statistical inference. If
$(X_{j})$ is a Gaussian process, then the Fourier transform
\[
S_{n}(\theta)=\sum_{j=1}^{n}X_{j}\exp(ji\theta)
\]
is complex Gaussian. \citet{fisher29} proposed a test for hidden
periodicities and obtained a distributional theory based on i.i.d.
Gaussian random variables. If $(X_{j})$ is not Gaussian, the
distribution of $S_{n}(\theta)$ typically does not have a close form
and one needs to resort to asymptotics. It is well known since
\citet{WW41} [see also \citet{LT08}] that for any stationary
sequence $(X_{j})_{j\in\mathbb{Z}}$ in $\mathcal{L}^{1}$ (namely
$E|X_{0}|<\infty$) there is a set $\Omega^{\prime}$ of probability $1$
such that for all $\theta$ and $\omega\in\Omega^{\prime}$,
$S_{n}(\theta)/n$ converges. Our problem is to investigate the speed of
this convergence by providing a central limit theorem for the real and
imaginary parts of ${S_{n}(\theta)}/\sqrt{n}$.

The above central limit problem was considered by many authors under
various dependence conditions. We mention Rosenblatt
[(\citeyear{Rosenblatt85}), Theorem 5.3, pa\-ge~131] who considered
mixing processes; Brockwell and Davis [(\citeyear{BD91}),
Theorem~10.3.2, page 347], \citet{walker65} and \citet{TH94}
discussed linear processes, and \citet{wu05} treated mixingales.
Other contributions can be found in \citet{Olshen67},
\citet{R76}, \citet{Yajima89}, \citet{Woodroofe92},
\citet{walker00}, \citet{Lahiri03} and \citet{LiuLin09}
among others.

To establish an asymptotic theory for ${S_{n}(\theta)}/\sqrt{n}$,
we shall provide the framework of stationary processes that can be
introduced in several equivalent ways. We assume that
$(\xi_{n})_{n\in\mathbb{Z}}$ is a stationary ergodic Markov chain
defined on a \mbox{probability} space $(\Omega,\mathcal{F},P)$ with
values in a measurable space. The marginal distribution is denoted
by $\pi(A)=P(\xi_{0}\in A)$. Next let $\mathcal{L}_{0}^{2}(\pi) $
be the set of functions such that $\int h^{2}\,d\pi<\infty$ and
$\int h\,d\pi=0$. Denote by $\mathcal{F}_{k}$ the $\sigma$-field
generated by $\xi_{j}$ with $j \leq k$, $X_{j} =h(\xi_{j})$. For
any integrable random variable $X$ we denote $E_{k}(X) =
E(X|\mathcal{F}_{k})$. We assume $h\in\mathcal{L}_{0}^{2}(\pi);$
in other words we assume $\| X_{0} \| := (E|X_{0}^{2}|)^{1/2} <
\infty$ and $E (X_{0}) = 0$. Notice that any stationary sequence
$(Y_{k})_{k\in\mathbb{Z}}$ can be viewed as a function of a Markov
process $\xi_{k} = (Y_{j}; j\leq k)$ with the function $g(\xi_{k})
= Y_{k}$.

The stationary stochastic processes may be also introduced in the
following alternative way. Let $T\dvtx\Omega\mapsto\Omega$ be a
bijective bi-measurable transformation preserving the probability.
Let $\mathcal{F}_{0}$ be a $\sigma$-algebra of $\mathcal{F}$
satisfying $\mathcal{F}_{1}\subseteq T^{-1} (\mathcal{F}_{0})$. We
then define the nondecreasing filtration $( \mathcal{F}_{j })_{j
\in\mathbb{Z}}$ by $\mathcal{F}_{j}=T^{-j}(\mathcal{F} _{0})$
(referred to as the stationary filtration). Let $X_{0}$ be a
random variable which is $\mathcal{F}_{0}$-measurable. We also
define the stationary sequence $(X_{j})_{j \in\mathbb{Z}}$ by
$X_{j} = X_{0} \circ T^{j}$. In this paper we shall use both
frameworks.

The rest of the paper is structured as follows. The main results
are presented in Section \ref{sec:main} and proved in Section
\ref{sec:proof}. Our proofs in Section \ref{sec:proof} provide an
interesting blend of harmonic analysis, martingale approximation
and theory of stationary processes. Examples of regular processes
and further extensions are given in Section \ref{sec:ex}.

\section{Main results}
\label{sec:main}

We shall assume that the following regularity condition holds:
%
\begin{equation}\label{regular}
E(X_{0}|\mathcal{F}_{-\infty})
=0,\qquad P\mbox{-almost surely},
\end{equation}
and also that the sequence is stationary and ergodic. The
regularity condition is quite mild and it is satisfied for many
popular processes used in practice. Section~\ref{sec:ex} provides
examples of stationary ergodic processes for which (\ref{regular})
holds.

We shall present first a central limit theorem for almost all
frequencies. In Theorem \ref{SD1}, we let the parameter $\theta$
be in the space $[0,2\pi]$, endowed with Borelian sigma algebra
and Lebesgue measure $\lambda$. We denote by ``$\Rightarrow$'' the
weak convergence, or convergence in distribution.
\begin{theorem}
\label{SD1} Let $(X_{k})_{k\in\mathbb{Z}}$ be a stationary ergodic
process such that (\ref{regular}) is satisfied. Then for almost
all $\theta\in(0,2\pi)$, the following convergence holds:
%
\begin{equation}\label{def1}%
\lim_{n\rightarrow\infty}{\frac{{E|S_{n}(\theta
)|^{2}}}{n}}=g(\theta
)\qquad\mbox{(say),}
\end{equation}
where $g$ is integrable over $\theta\in\lbrack0,2\pi]$, and
%
\begin{equation}\label{eq:asclt}%
\frac{1}{\sqrt{n}}[\mathrm{Re}(S_{n}(\theta)), \mathrm{Im}(S_{n}%
(\theta))]\Rightarrow\lbrack N_{1}(\theta), N_{2}(\theta)]\qquad
\mbox{under }P,
\end{equation}
where $N_{1}(\theta)$ and $N_{2}(\theta)$ are independent
identically distributed normal random variables mean $0$ and
variance $g(\theta)/2$.
\end{theorem}

As implied by Lemma \ref{key2} in Section \ref{sec:proof},
$g(\theta) / (2 \pi)$ is actually the spectral density associated
with the spectral distribution function induced by the covariances
%
\begin{equation}\label{cov}%
c_{j}=\mathrm{cov}(X_{0}, X_{j}),\qquad j\in\mathbb{Z}.
\end{equation}
More specifically, by Herglotz's theorem [\citet{BD91}], there exists a nondecreasing function $G$ (the spectral
distribution function) on $[0, 2\pi]$ such that, for all $j
\in\mathbb{Z}$,
%
\begin{equation}
\label{eq:June6}
c_{j} = \int_{0}^{2 \pi} \exp(i j \theta) \,d G(\theta).
\end{equation}
Hence, by Lemma \ref{key2}, $G$ is absolutely continuous and the
spectral density $G^{\prime}(\theta)$ equals to $g(\theta) / (2
\pi)$ almost surely. By (\ref{eq:June6}) or (\ref{eq:May312}),
$\int_{0}^{2\pi} g(\theta) \,d \theta= 2 \pi c_{0}$. So a nice
implication of our results is that, under the regularity
condition, we obtain an interesting representation of the spectral
densities (see Lemma \ref{key2} for details).

Following the proof of Theorem \ref{SD1}, by the Cram\'{e}r--Wold
device, for
\[
V_{n}(\omega,\theta)=\frac{S_{n}}{\sqrt{n}}(\omega,\theta),
\]
we have that, for almost all pairs $(\theta^{\prime},
\theta^{\prime\prime})$ (Lebesgue), $V_{n}(\omega,
\theta^{\prime})$ and $V_{n}(\omega, \theta^{\prime\prime})$ are
asymptotically independent. In this sense Theorem \ref{SD1}
justifies the folklore in the spectral domain analysis of time
series: the Fourier transforms of stationary processes are
asymptotically independent Gaussian. Namely, in the spectral or
Fourier domain, the Fourier-transformed processes are
asymptotically independent, while the original process can be very
strongly dependent (see Example \ref{ex:lp}).

Theorem \ref{SD1} substantially improves the result in \citet{wu05}
that proves (\ref{eq:asclt}) under the following stronger
condition:
%
\begin{equation}\label{eq:wucond}%
\sum_{n=1}^{\infty}{\frac{{\| E(X_{n}|\mathcal{F}_{0})\|
^{2}}}{n}} <
\infty.
\end{equation}

We shall also establish the following characterization of the
annealed CLT. Let $\operatorname{Id}_{2}$ denote the identity $2\times2$ matrix.
\begin{theorem}
\label{annealed} Under the same conditions as in Theorem \ref{SD1}
on the product space $([0,2\pi]\times\Omega, \mathcal{B}
\times\mathcal{F}, \lambda\times P)$ we have
%
\begin{equation}\label{annealed1}\qquad
\frac{1}{\sqrt{n}}[\mathrm{Re}(S_{n}(\theta)), \mathrm{Im}(S_{n}
(\theta))]\Rightarrow\lbrack g(U)/2]^{1/2}N(0,\operatorname{Id}_{2})\qquad
\mbox{under } \lambda\times P.
\end{equation}
Here $U$ is a random variable independent of $N(0,\operatorname{Id}_{2})$ and
uniformly distributed on $[0,2\pi]$ and $g(\cdot)$ is defined by
(\ref{def1}).
\end{theorem}

Two types of stochastic processes can be considered concerning the
partial sum $S_{n}(\theta)$. The process $V_{n}(\omega,\theta)$
indexed by $\theta$ is asymptotically Gaussian white noise. For
another version, we consider
\[
W_{n}(t,\omega,\theta) = \frac{S_{ \lfloor n t \rfloor}}{\sqrt{n}}
(\omega,\theta),\qquad 0 \le t \le1,
\]
where $\lfloor x \rfloor= \max\{ k \in\mathbb{Z}\dvtx k \le x \}$
is the integer part of $x$. We shall prove the following
invariance principle:
\begin{proposition}
\label{IP} Assume that $(X_{k})$ is stationary ergodic and
satisfies (\ref{regular}). Then $W_{n}(t,\omega,\theta)$ is tight
in $D(0,1)$, and
\begin{eqnarray*}
&&\lbrack\mathrm{Re}(W_{n}(t,\omega,\theta)),
\mathrm{Im}(W_{n}(t,\omega,\theta))]\\
&&\qquad\Rightarrow\lbrack
g(U)/2]^{1/2}[W^{\prime}(t),W^{\prime\prime}(t)]\qquad\mbox{under
}\lambda\times P,
\end{eqnarray*}
where $(W^{\prime}(t),W^{\prime\prime}(t))$ are two independent
standard Brownian motions independent of $U$, and $U$ is a random
variable uniformly distributed on $[0,2\pi]$.
\end{proposition}

We now give some remarks and discussions.
\begin{remark}[(Nonadapted case)]
Our CLT also holds if $X_{0}$ is not
$\mathcal{F}_{0}$-measurable. Then clearly the regularity
condition we shall impose is
\[
X_{0}\mbox{ is
}\mathcal{F}_{\infty}\mbox{-measurable}\quad\mbox{and}\quad
E(X_{0}|\mathcal{F}_{-\infty})=0\qquad\mbox{almost surely}.
\]
\end{remark}
\begin{remark}[(Adapted nonregular case)]
For general adapted sequences our
CLT result still holds under centering. Let $\tilde{S}_{n}(\theta)
= S_{n} (\theta)-E(S_{n}(\theta) |\mathcal{F}_{-\infty})$ and
$\tilde{X}_{k} = X_{k}-E(X_{k}|\mathcal{F}_{-\infty})$, where
$E(X_{k}|\mathcal{F}_{-\infty})$ denotes the following limit that
holds almost surely and in $\mathcal{L}^{2}$
\[
\lim_{n\rightarrow\infty}E(X_{k}|\mathcal{F}_{-n})=E(X_{k}|\mathcal
{F}%
_{-\infty}).
\]
Then, for almost all $\theta\in(0,2\pi)$, we have
\[
\frac{1}{\sqrt{n}}[\mathrm{Re}(\tilde{S}_{n}(\theta)), \mathrm
{Im}(\tilde
{S}_{n}(\theta))]\Rightarrow
N\biggl(0,{\frac{{g(\theta)}}{2}}\operatorname{Id}_{2}\biggr)\qquad\mbox{under }P
\]
by applying Theorem \ref{SD1} to the stationary sequence
$\tilde{X}_{k}$. Therefore the conclusion of Theorem \ref{SD1}
holds if we replace the assumption of regularity (\ref{regular})
by the following: for $\lambda$-almost all $\theta$
\[
\frac{1}{\sqrt{n}}E(S_{n}(\theta)|\mathcal{F}_{-\infty})
\to0\qquad\mbox{in probability}.
\]
Now, since $\|E(S_{n}(\theta)|\mathcal{F}_{-\infty})\|_{2} \leq
\|E(S_{n}(\theta)|\mathcal{F}_{-n}) \|_{2} \leq\|
E(S_{n}(\theta)|\mathcal{F}_{0})\|_{2}$, Theorem~\ref{SD1} still
holds if
\[
\frac{1}{\sqrt{n}}\|E(S_{n}(\theta)|\mathcal{F}_{-n})\|_{2}
\rightarrow0\qquad\mbox{as }n\rightarrow\infty
\]
or under the condition
%
\begin{equation}\label{condfn}
\frac{1}{\sqrt{n}}\|E(S_{n}(\theta)|\mathcal{F}_{0})\|_2
\rightarrow0 \qquad\mbox{as }n\to\infty.
\end{equation}
\end{remark}
\begin{remark}[(Conditional CLT)]
Since we use in the proof martingale
approximation actually our CLT is a conditional CLT, that allows
for a random change of measure. See \citet{HH80} and
\citet{DM02}.
\end{remark}
\begin{remark}[(Periodogram)]
We notice that, as a consequence of Theorem~\ref{SD1},
for sequences satisfying (\ref{regular}) the
periodogram $n^{-1}|S_{n} (\theta)|^{2}$ is asymptotically
distributed as ${\frac{{g(\theta)}}{2}} \chi^{2}(2)$ for almost
all frequencies.
\end{remark}
\begin{remark}[(Resulting identities)]
By using the invariance principle in
Proposition \ref{IP} we can get the convergence of many
interesting functionals of $S_{n}(\theta)$ and periodograms. As a
consequence of Proposition \ref{IP} we can get, for instance,
\begin{eqnarray*}
&&\frac{1}{\Vert X_{0}\Vert^{2}n\pi}\int_{0}^{2\pi}E \Biggl[
\max_{1\leq m\leq n}\sum_{k=1}^{m}X_{k}\cos(k\theta) \Biggr]
^{2}\,d\theta\\
&&\qquad \rightarrow
E \Bigl\vert\sup_{0\leq t\leq1}W(t) \Bigr\vert^{2}
=\int_{0}^{\infty}2\bigl[1-\Phi\bigl(\sqrt{y}\bigr)\bigr]\,dy=1
\end{eqnarray*}
by noting that $P(\sup_{0\leq t\leq1}W(t)\ge u) = 2 P(W(1)\ge u)$
for $u\ge0$. Here $\Phi(\cdot)$ is the standard Gaussian
distribution function.
\end{remark}

\section{Examples}
\label{sec:ex} Here we present several examples of processes for
which the conclusions of Theorems \ref{SD1} and \ref{annealed}
hold.

Clearly condition (\ref{regular}) is satisfied if the left tail
sigma field $\mathcal{F}_{-\infty}$ is trivial. These processes
are called regular [see Chapter 2, Volume 1 in \citet{Bradley07}].
Notice, however, that our condition (\ref{regular}) refers rather to
the function $X_{0}=f(\xi_{0})$ in relation to the tail field
$\mathcal{F}_{\infty}$.
\begin{example}[(Mixing sequences)]
We shall introduce the following mixing
coefficients: for any two $\sigma$-algebras $\mathcal{A}$ and
$\mathcal{B}$ define the strong mixing coefficient
\[
\alpha(\mathcal{A},\mathcal{B)=}\sup\{|P(A\cap B)-P(A)P(B)|\dvtx
A\in\mathcal{A}, B\in\mathcal{B\}}
\]
and the $\rho$-mixing coefficient, also known as maximal
coefficient of correlation
\[
\rho(\mathcal{A},\mathcal{B})=\sup\{\mathrm{Cov}(X,Y)/\Vert
X\Vert_{2}\Vert Y\Vert_{2}\dvtx
X\in\mathcal{L}^{2}(\mathcal{A}), Y\in\mathcal{L}^{2}
(\mathcal{B})\}.
\]
For the stationary sequence of random variables
$(X_{k})_{k\in\mathbb{Z}}$, $\mathcal{F}^{n}$ denotes the
$\sigma$-field generated by $X_{i}$ with indices $i\geq n$, and
$\mathcal{F}_{m}$ denotes the $\sigma$-field generated by $X_{i}$
with indices $i\leq m$. The sequences of coefficients $\alpha(n)$
and $\rho(n)$ are then defined by
\[
\alpha(n)=\alpha(\mathcal{F}_{0},\mathcal{F}^{n})
\quad\mbox{and}\quad
\rho(n)=\rho(\mathcal{F}_{0},\mathcal{F}^{n}),
\]
respectively. For strongly mixing sequences, namely the strong
mixing coefficients $\alpha(n)\rightarrow0$, the tail sigma field
is trivial [see Claim 2.17a in \citet{Bradley07}]. Examples of this
type include Harris recurrent Markov chains. If
$\lim_{n\rightarrow\infty}\rho(n)<1$, then the tail sigma field is
also trivial [see Proposition 5.6 in \citet{Bradley07}].
\end{example}
\begin{example}[(Functions of Gaussian processes)]
Assume $(Y_{k})$ is a
stationary Gaussian sequence and define $X_{n}=f(Y_{k}, k \leq
n)$. Let $f$ be such that $E(X_{0})=0$ and $E(X_{0}^{2})<\infty$.
Since any Gaussian sequence can be represented as a function of
i.i.d. random variables, the process is then regular.
\citet{Rosenblatt81} considered Fourier transforms of functionals of Gaussian
sequences.
\end{example}
\begin{example}[(Functions of i.i.d. random variables)]
\label{ex:lp}
Let $\varepsilon_k $ be i.i.d. and consider $X_n = f(\varepsilon_k, k
\le n)$. These are regular processes and therefore Theorems~\ref{SD1}
and \ref{annealed} are applicable. Examples include
linear processes, functions of linear processes and iterated
random functions [\citet{WW00}, among others]. For
example, let $X_{n} = \sum_{j=0}^{\infty} a_{j} \varepsilon_{n-j}$,
where $\varepsilon_{j}$ are i.i.d. with mean $0$ and variance $1$,
and $a_{j}$ are real coefficients with
$\sum_{j=1}^{\infty}a_{j}^{2} < \infty$. In this case $X_{n}$ is
well defined, and, by Lemmas \ref{key} and \ref{key2}, the spectral
density is $g(\theta) / (2 \pi)$, where
\[
g(\theta) = \Biggl| \sum_{j=0}^{\infty}a_{j} \exp(i j \theta)
\Biggr| ^{2}.
\]
As a specific example, let $a_{j} = j^{-1/2} / \log j$, $j \ge2$,
and $a_{0} = a_{1} = 1$. By elementary manipulations, the
covariances $c_{j} \sim(\log j)^{-1}, $ which decays very slowly
as $j \to\infty$, hence suggesting strong dependence. For this
example, condition (\ref{eq:wucond}) is violated. By the Tauberian
theorem, as $\theta\to0$, $g(\theta) \sim\pi/ ({|\theta| \log^{2}}
|\theta| ), $ which has a pole at $\theta= 0$.
\end{example}
\begin{example}[(Reversible Markov chains)]
As before, let $\xi_{j}$ be a
stationary ergodic Markov chain with values in a measurable space.
We use the notation and constructions from the \hyperref[sec1]{Introduction}. The
marginal distribution and the transition metric are denoted by
$\pi(A)=P(\xi_{0}\in A)$ and $Q(\xi_{0},A)=P(\xi_{1}\in
A|\xi_{0})$. In addition $Q$ denotes the operator $Qf(\xi)=\int
f(z)Q(\xi,dz)$. Let $Q^{\ast}$ be the adjoint operator of the
restriction of $Q$ to $\mathcal{L}^{2}(\pi)$ and assume
$Q=Q^{\ast}$. Then, for any $f \in\mathcal{L}_0^2(\pi)$ the
central limit theorem of Theorem \ref{SD1} holds. To see this we
shall verify condition (\ref{condfn}). By spectral calculus
\[
\Vert
E(S_{n}(\theta)|\mathcal{F}_{0})\Vert^{2}
=\int_{-1}^{1}\Biggl|\sum_{k=1}^{n}(t\exp(i\theta))^{k}\Biggr|^{2}
\rho_{f}(dt),
\]
where $\rho_{f}$ denotes the spectral measure of $f$ with respect
to $Q$ [see, e.g., \citet{BI94} for this
identity]. For $\theta\neq0$, $\pi$ and $-1\leq t\leq1$, we have
\begin{eqnarray*}
\Biggl|\sum_{k=1}^{n}(t\exp(i\theta))^{k}\Biggr|^{2}
&\le& 4|1-\exp(i\theta)t|^{-2} \\
&=&4(1+t^{2}-2t\cos\theta)^{-1}
\le4\bigl(1-(\cos\theta)^{2}\bigr)^{-1}.
\end{eqnarray*}
Therefore, for $\lambda$-almost all $\theta$%
\[
\frac{1}{n}\|E(S_{n}(\theta)|\mathcal{F}_{0})\|^{2}
\rightarrow0.
\]
\end{example}

\section{Proofs}
\label{sec:proof} We shall establish first some preparatory
lemmas. The almost sure convergence in Lemma \ref{key} was shown
by \citet{wu05}. The convergence in $\mathcal{L}^{2}$ is new here.
For $k\in\mathbb{Z}$ we define the projection operator by
%
\begin{equation}\label{proj}%
\mathcal{P}_{k}\cdot=E(\cdot|\mathcal{F}_{k}) -E(\cdot|\mathcal{F}
_{k-1}).
\end{equation}
\begin{lemma}
\label{key} Let
\[
T_{n}(\theta)=\sum_{j=0}^{n}X_{j}\exp(ji\theta)=X_{0}+S_{n}(\theta).
\]
Under (\ref{regular}), for $\lambda$-almost all $\theta$
(Lebesgue), we have
%
\begin{equation}\label{DEF1}
\mathcal{P}_{0}T_{n}(\theta)\rightarrow\sum_{l=0}^{\infty}\mathcal
{P}_{0}%
X_{l}\exp(li\theta)=:D_{0}(\theta),\qquad P\mbox{-almost
surely and in
}\mathcal{L}^{2}.\hspace*{-32pt}
\end{equation}
\end{lemma}
\begin{pf}
By (\ref{regular}), $\sum_{k\in\mathbb{Z}} \Vert
\mathcal{P}_{k} X_0 \Vert_{2}^{2} = \Vert X_{0}\Vert_{2}^{2} <
\infty$, we have $\sum_{k \in\mathbb{Z}} |\mathcal{P}_0 X_k|^2 <
\infty$, $P$-almost surely. Therefore by Carleson's (\citeyear{Carleson66}) theorem,
for almost all $\omega$, $\sum_{1\leq k\leq n}(\mathcal{P}_{0}
X_{k})\exp(i k \theta)$ converges $\lambda$-almost surely, where
$\lambda$ is the Lebesgue measure on $[0,2\pi]$. Denote the limit
by $D_{0}=D_{0}(\theta)$. We now consider the set
\[
A=\bigl\{(\theta,\omega)\subset\lbrack0,2\pi]\times\Omega,\mbox{ where
}\{\mathcal{P}_{0}S_{n}(\theta)\}_{n}\mbox{ does not converge}\bigr\}
\]
and notice that almost all sections for $\omega$ fixed have
Lebesgue measure $0$. So by Fubini's theorem the set $A$ has
measure $0$ in the product space and therefore, again by Fubini's
theorem, almost all sections for $\theta$ fixed have probability
$0$. It follows that for almost all $\theta$,
$\mathcal{P}_{0}(S_{n}(\theta))\rightarrow D_{0}$ almost surely
under $P$. Next, by the maximal inequality in \citet{HY74},
there is a constant $C$ such that
\[
\int_{0}^{2\pi} \Bigl[
{\sup_{n}}|\mathcal{P}_{0}(S_{n}(\theta))|^{2} \Bigr]
\lambda(d\theta)\leq C\sum_{k}|\mathcal{P}_{0}X_{k}|_{2}^{2},%
\]
and then we integrate
\[
\int_{0}^{2\pi}E \Bigl[
{\sup_{n}}|\mathcal{P}_{0}(S_{n}(\theta))|^{2} \Bigr]
\lambda(d\theta)\leq C\Vert X_{0}\Vert_{2}^{2}<\infty.
\]
Therefore,
\[
E \Bigl[ {\sup_{n}}|\mathcal{P}_{0}(S_{n}(\theta))|^{2} \Bigr]
<\infty\qquad\mbox{for almost all }\theta.
\]
Since
$|\mathcal{P}_{0}S_{n}(\theta)|<\sup_{n}|\mathcal
{P}_{0}S_{n}(\theta)|$,
and the last one is integrable for almost all $\theta$, by the
Lebesgue dominated convergence we have that
$\mathcal{P}_{0}(S_{n}(\theta))$ converges in $\mathcal{L}^{2}$.
\end{pf}
\begin{lemma}
\label{key2} Let $g(\theta)=E|D_{0}(\theta)|^{2}$. For all
$j\in\mathbb{Z} $, we have
%
\begin{equation}\label{eq:May312}%
\int_{0}^{2\pi}g(\theta)\exp(ji\theta)\,d\theta=2\pi c_{j},
\end{equation}
where $c_{j}$ are defined by (\ref{cov}). So $(c_{j})$ are the
Fourier coefficients of $g$. Additionally, for almost all
$\theta$,
%
\begin{equation}\label{may311}
\lim_{n\to\infty}\frac{E|S_{n}(\theta)|^{2}}{n}=g(\theta).
\end{equation}
\end{lemma}
\begin{pf}
Without loss of generality we let $j\geq0$ and $j<n$.
As before, let
\[
\mathcal{P}_{0}T_{n}(\theta)=\sum_{l=0}^{n}\mathcal{P}_{0}X_{l}\exp
(li\theta).
\]
By elementary trigonometric identities, we have
%
\begin{equation}\label{eq:May313}%
{\frac{1}{{2\pi}}}\int_{0}^{2\pi}|\mathcal{P}_{0}T_{n}(\theta
)|^{2}%
\exp(ji\theta)\,d\theta=\sum_{l=0}^{n-j}(\mathcal
{P}_{0}X_{l})(\mathcal{P}%
_{0}X_{l+j}).
\end{equation}
Since $X_{j}=\sum_{l\in\mathbb{Z}}\mathcal{P}_{l}X_{j}$, by
orthogonality of martingale differences and stationarity we have
that
\begin{eqnarray*}
c_{j} &=&\lim_{N \to\infty} E \Biggl[ \Biggl( \sum_{l=-N}^{0}
\mathcal{P}_{l}X_{0} \Biggr) \Biggl(
\sum_{l=-N}^{0}\mathcal{P}_{l} X_{j} \Biggr) \Biggr] \\
&=&\lim_{N\to\infty}\sum_{l=-N}^{0}
E[(\mathcal{P}_{l}X_{0})(\mathcal{P}_{l}X_{j})]\\
&=&\lim_{n \to\infty}\sum_{l=0}^{n-j}E[(\mathcal{P}_{0}X_{l}
)(\mathcal{P}_{0}X_{l+j})].
\end{eqnarray*}
By (\ref{eq:May313}) and the Lebesgue dominated convergence
theorem, as in the proof of Lemma \ref{key}, (\ref{eq:May312})
follows in view of Hunt's maximal inequality since
${\sup_{n}}|\mathcal{P}_{0}\times T_{n}(\theta)|$ is integrable.

Now we prove (\ref{may311}). By stationarity, we have
%
\begin{eqnarray}
\label{May315}
\frac{1}{n}E|S_{n}(\theta)|^{2}
&=&\frac{1}{n}\sum_{j=1}^{n}\sum_{l=1}^{n}
E(X_{j}X_{l})\exp(ij\theta)\exp(-il\theta)\nonumber\\
&=&\frac{1}{n}\sum_{j=1}^{n-1}
\sum_{l=1}^{n}c_{j-l}\exp\bigl((j-l)i\theta\bigr)\\
&=&\sum_{j=-(n-1)}^{n-1}\biggl(1-\frac{|j|}{n}\biggr)c_{j}\exp(ij\theta).\nonumber
\end{eqnarray}
Namely $E|S_{n}(\theta)|^{2}/n$ is the Cesaro average of the sum
$\sum_{j=-l}^{l}c_{j}\exp(ij\theta)$. Note that $g(\theta)=\Vert
D_{0}(\theta)\Vert^{2}$ is integrable over $[0,2\pi]$. Therefore
by the Fej\'{e}r--Lebesgue theorem [cf. \citet{Bary64}, page 139 or
Theorem 15.7 in Champeney (\citeyear{Champeney89})], (\ref{may311}) holds for
$\lambda$-almost all $\theta\in\lbrack0,2\pi]$ (Lebesgue).
\end{pf}
\begin{remark}
In the proof of Lemma \ref{key2}, (\ref{May315}) implies that
the sequence $(c_{j}\exp(ij\theta))$ is Cesaro summable. It turns
out that, generally speaking, $\sum_{j=0}^{\infty} c_{j}\times \exp(i j
\theta)$ may not exist for almost all $\theta$. Consider the
example in \citet{Kolmogorov23} [see also Theorem 3.1, page 305, in
\citet{Zygmund02}]: there exists a sequence of nonnegative
trigonometric polynomials $f_{n}$ with constant term $1/2$,
a~sequence of positive integers $q_{k}\rightarrow\infty$ and a
positive sequence $A_{n}\rightarrow\infty$, such that the function
\[
g(x)=\sum_{k=1}^{\infty}{\frac
{{f_{n_{k}}(q_{k}x)}}{{A_{n_{k}}^{1/2}}}}%
\]
is integrable. However, for almost all $\theta$, the Fourier sum
$\sum_{l=1}^{\infty}c_{l}\exp(li\theta)$ diverges, where
$c_{l}=\int_{0}^{2\pi}g(\theta)\exp(li\theta)\,d\theta$ is the
Fourier coefficient of $g$. Let $G(x)=\int_{0}^{x}g(u)\,du$. By
Herglotz's theorem [\citet{BD91}], there exists a
stationary process $(X_{j})$ such that its spectral distribution
function is $G$ and its covariance function is $c_{l}$.
\end{remark}
\begin{lemma}
\label{UI} Assume (\ref{regular}). On the product space
$([0,2\pi]\times\Omega, \mathcal{B}\times\mathcal{F},
\lambda\times P)$ we have that
\[
\biggl(\frac{{\max_{1\leq k\leq
n}}|S_{k}(\theta)|^{2}}{n}\biggr)_{n\geq1}\qquad
\mbox{is uniformly integrable.}%
\]
\end{lemma}
\begin{pf}
Let $m$ be a positive integer. We shall decompose
the partial sums in a sum of $m$ martingales and a remainder in
the following way:
\[
\frac{S_{k}(\theta)}{\sqrt{n}}=\frac{1}{\sqrt{n}}\sum
_{l=1}^{k}\exp
(il\theta)\sum_{j=0}^{m-1}\mathcal{P}_{l-j}(X_{l})+\frac{1}{\sqrt
{n}}%
\sum_{l=1}^{k}\exp(il\theta)E_{l-m}(X_{l}).
\]
Notice that for any $0\leq j\leq m-1$,%
\[
\sum_{l=1}^{k}\exp(il\theta)\mathcal{P}_{l-j}(X_{l})
\]
is a martingale adapted to the filtration $(\mathcal{B}\times\mathcal
{F}%
_{k})$. Moreover, since $(X_{k})_{k\in\mathbb{Z}}$ is a stationary
sequence
with variables square integrable, it follows that $(X_{k}^{2})_{k\in
\mathbb{Z}%
}$ is a uniformly integrable sequence. This fact implies that for
$j$ fixed the sequence
$(\mathcal{P}_{k-j}(X_{k}))_{k\in\mathbb{Z}}$ is also uniformly
integrable. It follows that $\sum_{l=1}^{k}\exp(il\theta)\times\mathcal{P}
_{l-j}(X_{l})$ is a martingale with uniformly integrable
differences under the measure $\lambda\times P$. It is known that
for a martingale with uniformly integrable differences we have
\[
\frac{1}{n}\max_{1\leq k\leq n}\Biggl|\sum_{l=1}^{k}\exp(il\theta
)\mathcal{P}%
_{l-j}(X_{l})\Biggr|^{2}%
\]
is uniformly integrable [see, e.g., \citet{DR00},
Proposition 1].

The result follows since by \citet{HY74} maximal
inequality
\[
\int_{0}^{2\pi}\max_{1\leq k\leq n}\Biggl|\sum_{l=1}^{k}\exp(il\theta
)E_{l-m}%
(X_{l})\Biggr|^{2}\,d\theta\leq {C\sum_{l=1}^{n}}|E_{l-m}(X_{l})|^{2},%
\]
and therefore, denoting by $\mathbf{E}$ the expected value with
respect to $\lambda\times P$, we have by regularity condition
(\ref{regular}) that
\begin{eqnarray*}
\frac{1}{n}\mathbf{E} \Biggl[ \max_{1\leq k\leq
n}\Biggl|\sum_{l=1}^{k}\exp
(il\theta)E_{l-m}(X_{l})\Biggr|^{2} \Biggr] &\leq& \frac{C}{n}\sum_{l=1}%
^{n}E|E_{l-m}(X_{l})|^{2}\\
&=& CE|E_{-m}(X_{0})|^{2}\rightarrow0
\end{eqnarray*}
as $m\rightarrow\infty$ uniformly in $n$. So, the uniform
integrability follows.
\end{pf}

\subsection[Proof of Theorem 2.1]{Proof of Theorem \protect\ref{SD1}}

The first assertion of Theorem \ref{SD1} is just Lem\-ma~\ref{key}.
We now prove (\ref{eq:asclt}).

\textit{Step} 1. \textit{The construction of martingale.}

Define the projector operator by (\ref{proj}). Then we construct
as in Lemma \ref{key}
\[
\mathcal{P}_{1}(S_{n}(\theta))=E(S_{n}(\theta)|\mathcal
{F}_{1})-E(S_{n}%
(\theta)| \mathcal{F}_{0})=\sum_{k=1}^{n}\exp(i k \theta)\mathcal
{P}_{1}%
(X_{k})
\]
and then by Lemma \ref{key} for almost all $\theta$
\[
\mathcal{P}_{1}(S_{n}(\theta))\rightarrow D_{1}(\theta) \qquad\mbox{in
} \mathcal{L}^{2}.
\]
To verify it is a martingale we start from
\[
E(\mathcal{P}_{1}(S_{n}(\theta))|\mathcal{F}_{0})=0\qquad\mbox{almost
surely under }P
\]
and by the contractive property of the conditional expectation
\[
0=E(\mathcal{P}_{1}(S_{n}(\theta))|\mathcal{F}_{0})\rightarrow
E( D_{1}(\theta)|\mathcal{F}_{0}) \qquad\mbox{in }
\mathcal{L}^{2}.
\]
We then construct the sequence of stationary martingale
differences $(D_{k}(\theta))_{k \ge1}$, given by
\[
\mathcal{P}_{k}\bigl(S_{n+k}(\theta)-S_{k}(\theta)\bigr)\rightarrow
\exp (i k \theta)D_{k}(\theta) \qquad\mbox{in } \mathcal{L}^{2}.
\]

\textit{Step} 2. \textit{Martingale approximation.}

Denote by
\[
M_{n}(\theta)=\sum_{1\leq k\leq n}\exp(ik\theta)D_{k}(\theta).
\]
We show that, for almost all $\theta$,
%
\begin{equation}\label{eq:June12}%
\frac{E|S_{n}(\theta)-M_{n}(\theta)|^{2}}{n}\rightarrow0.
\end{equation}
To this end, note that
$S_{n}(\theta)-E(S_{n}(\theta)|\mathcal{F}_{0})$ and
$E(S_{n}(\theta)|\mathcal{F}_{0})$ are orthogonal, we have
%
\begin{equation}\label{eq:June11}%
\Vert S_{n}(\theta)\Vert^{2}=\Vert S_{n}(\theta)-E(S_{n}(\theta
)|\mathcal{F}%
_{0})\Vert^{2}+\Vert E(S_{n}(\theta)|\mathcal{F}_{0})\Vert^{2}.
\end{equation}
For those $\theta$ such that (\ref{DEF1}) holds, we have, by the
orthogonality
of martingale differences and the stationarity, that%
%
\begin{eqnarray}\label{June13}
&&\Vert
S_{n}(\theta)-E(S_{n}(\theta)|\mathcal{F}_{0})-M_{n}(\theta)\Vert
^{2}\nonumber\\
&&\qquad=\sum_{k=1}^{n}\bigl\Vert\mathcal{P}_{k}\bigl(S_{n}(\theta)-M_{n}(\theta
)\bigr)\bigr\Vert
^{2}
=\sum_{k=1}^{n}\Vert\mathcal{P}_{k}S_{n}(\theta)-e^{ik\theta}D_{k}%
(\theta)\Vert^{2}\\
&&\qquad=\sum_{k=1}^{n}\Vert\mathcal{P}_{0}T_{n-k}(\theta
)-D_{0}(\theta)\Vert^{2}=o(n).\nonumber
\end{eqnarray}
Hence, by (\ref{eq:June11}) and (\ref{may311}), we have
\[
\limsup_{n\rightarrow\infty}{\frac{{\Vert E(S_{n}(\theta)|\mathcal
{F}%
_{0})\Vert^{2}}}{n}}=\limsup_{n\rightarrow\infty}{\frac{{\Vert
S_{n}%
(\theta)\Vert^{2}-\Vert M_{n}(\theta)\Vert^{2}}}{n}}=0
\]
by noting that $\Vert M_{n}(\theta)\Vert^{2}=n\Vert
D_{0}(\theta)\Vert^{2}$. Hence we have (\ref{eq:June12}) in view
of (\ref{June13}).

\textit{Step} 3. \textit{The CLT for the approximating martingale}.

It remains just to prove central limit theorem for complex valued
martingale
\[
\frac{1}{\sqrt{n}}\sum_{1\leq k\leq n}\exp(ik\theta)D_{k}(\theta).
\]
As a matter of fact we shall provide a central limit theorem for
the real part and imaginary part and show that in the limit they
are independent. The proof was carefully written down in
\citet{wu05}. By the Cram\'{e}r--Wold device we have to study the limiting
distribution of the martingale
\[
\frac{1}{\sqrt{n}}\sum_{1\leq k\leq n}[s\operatorname{Re}(\exp(ik\theta
)D_{k}%
(\theta))+t\operatorname{Im}(\exp(ik\theta)D_{k}(\theta))].
\]
By the Raikov-type of argument, in order to prove the CLT we have
only to show
\[
\frac{1}{n}\sum_{1\leq k\leq n}[s\operatorname{Re}(\exp(ik\theta)D_{k}%
(\theta))+t\operatorname{Im}(\exp(ik\theta)D_{k}(\theta
))]^{2}\rightarrow^{p}%
\frac{(s^{2}+t^{2})\sigma^{2}(\theta)}{2}.
\]
This follows from combining the following two facts. First by
stationarity
\[
{\frac{1}{n}\sum_{1\leq k\leq n}}|\exp(ik\theta)D_{k}(\theta
)|^{2}={\frac{1}%
{n}\sum_{1\leq k\leq n}}|D_{k}(\theta)|^{2}\rightarrow E|D_{0}(\theta
)|^{2}%
\]
and then by Lemma 5 in \citet{wu05} for almost all $\theta$
\[
\frac{1}{n}\sum_{1\leq k\leq
n}[\exp(ik\theta)D_{k}(\theta)]^{2}\rightarrow0.
\]
The rest is simple algebra.

\subsection[Proof of Theorem 2.2]{Proof of Theorem \protect\ref{annealed}}

Consider the product space
$([0,2\pi]\times\Omega, \mathcal{B}\times
\mathcal{F}, \lambda\times P)$. Let $\mathbf{P}=\lambda\times P$
and $\mathbf{E}$ the corresponding expected value. We already have
shown in Theorem \ref{SD1}, that for $\lambda$-almost all
$\theta$,
\[
E\exp\biggl[ \frac{i}{\sqrt{n}}\bigl(s\operatorname{Re}(S_{n}(\theta))+t\operatorname{Im}
(S_{n}(\theta))\bigr) \biggr] \rightarrow\exp\biggl[
-\frac{(s^{2}+t^{2})g(\theta
)}{4} \biggr]
\]
for $s,t\in\mathbb{R}$. Then we integrate with $\theta$ and by the
dominated convergence theorem we obtain
\[
\mathbf{E}\exp\biggl[ \frac{i}{\sqrt{n}}\bigl(s\operatorname{Re}(S_{n}(\theta
))+t\operatorname{Im}(S_{n}(\theta))\bigr) \biggr]
\rightarrow{\frac{1}{{2\pi}}}\int_{0}^{2\pi}\exp\biggl[
-\frac{(s^{2}+t^{2})g(\theta)}{4} \biggr] \,d\theta.
\]
We then identify the limiting distribution as being a mixture of
two independent random variables: a standard normal variable with
a variable uniformly distributed on $[0,2\pi]$.

\subsection[Proof of Proposition 2.1]{Proof of Proposition \protect\ref{IP}}

It is easy to see that the finite-dimensional distributions are
convergent. So we just have to prove tightness. From
\citet{Billingsley99}, stationarity and standard considerations this follows by
Lemma \ref{UI}.

\section*{Acknowledgments}
The first author would like to thank
Sergey Utev for his conjecture related to Theorem \ref{SD1} that
motivated this paper and Richard Bradley for useful discussions.
The authors are grateful to the referee for carefully reading the
manuscript and useful suggestions.


%
\printaddresses


\begin{thebibliography}{34}

\bibitem[\protect\citeauthoryear{Bary}{1964}]{Bary64}
\begin{bbook}[mr]
\bauthor{\bsnm{Bary},~\bfnm{N.~K.}\binits{N.~K.}}
(\byear{1964}).
\btitle{A Treatise on Trigonometric Series}. 
\bpublisher{Macmillan}, \baddress{New York}.
\bid{mr={0171116}}
\end{bbook}
\endbibitem

\bibitem[\protect\citeauthoryear{Billingsley}{1999}]{Billingsley99}
\begin{bbook}[mr]
\bauthor{\bsnm{Billingsley},~\bfnm{Patrick}\binits{P.}}
(\byear{1999}).
\btitle{Convergence of Probability Measures},
\bedition{2nd} ed.
\bpublisher{Wiley}, \baddress{New York}.
\bid{doi={10.1002/9780470316962}, mr={1700749}}
\end{bbook}
\endbibitem

\bibitem[\protect\citeauthoryear{Bradley}{2007}]{Bradley07}
\begin{bbook}[mr]
\bauthor{\bsnm{Bradley},~\bfnm{Richard~C.}\binits{R.~C.}}
(\byear{2007}).
\btitle{Introduction to Strong Mixing Conditions} 
\bvolume{1, 2, 3}.
\bpublisher{Kendrick Press}, \baddress{Heber City, UT}.
\bid{mr={2325294}}
\end{bbook}
\endbibitem

\bibitem[\protect\citeauthoryear{Borodin and Ibragimov}{1994}]{BI94}
\begin{barticle}[mr]
\bauthor{\bsnm{Borodin},~\bfnm{A.~N.}\binits{A.~N.}} \AND
  \bauthor{\bsnm{Ibragimov},~\bfnm{I.~A.}\binits{I.~A.}}
(\byear{1994}).
\btitle{Limit theorems for functionals of random walks}.
\bjournal{Tr. Mat. Inst. Steklova}
\bvolume{195}
\bpages{286}.
\bid{mr={1368394}}
\end{barticle}
\endbibitem

\bibitem[\protect\citeauthoryear{Brockwell and Davis}{1991}]{BD91}
\begin{bbook}[mr]
\bauthor{\bsnm{Brockwell},~\bfnm{Peter~J.}\binits{P.~J.}} \AND
  \bauthor{\bsnm{Davis},~\bfnm{Richard~A.}\binits{R.~A.}}
(\byear{1991}).
\btitle{Time Series: Theory and Methods},
\bedition{2nd} ed.
\bpublisher{Springer}, \baddress{New York}.
\bid{doi={10.1007/978-1-4419-0320-4}, mr={1093459}}
\end{bbook}
\endbibitem

\bibitem[\protect\citeauthoryear{Carleson}{1966}]{Carleson66}
\begin{barticle}[mr]
\bauthor{\bsnm{Carleson},~\bfnm{Lennart}\binits{L.}}
(\byear{1966}).
\btitle{On convergence and growth of partial sumas of {F}ourier series}.
\bjournal{Acta Math.}
\bvolume{116}
\bpages{135--157}.
\bid{mr={0199631}}
\end{barticle}
\endbibitem

\bibitem[\protect\citeauthoryear{Champeney}{1989}]{Champeney89}
\begin{bbook}[mr]
\bauthor{\bsnm{Champeney},~\bfnm{D.~C.}\binits{D.~C.}}
(\byear{1989}).
\btitle{A Handbook of {F}ourier Theorems}.
\bpublisher{Cambridge Univ. Press}, \baddress{Cambridge}.
\bid{mr={900583}}
\end{bbook}
\endbibitem

\bibitem[\protect\citeauthoryear{Dedecker and Rio}{2000}]{DR00}
\begin{barticle}[mr]
\bauthor{\bsnm{Dedecker},~\bfnm{J{\'e}r{\^o}me}\binits{J.}} \AND
  \bauthor{\bsnm{Rio},~\bfnm{Emmanuel}\binits{E.}}
(\byear{2000}).
\btitle{On the functional central limit theorem for stationary processes}.
\bjournal{Ann. Inst. H. Poincar\'e Probab. Statist.}
\bvolume{36}
\bpages{1--34}.
\bid{doi={10.1016/S0246-0203(00)00111-4}, mr={1743095}}
\end{barticle}
\endbibitem

\bibitem[\protect\citeauthoryear{Dedecker and Merlev{\`e}de}{2002}]{DM02}
\begin{barticle}[mr]
\bauthor{\bsnm{Dedecker},~\bfnm{J{\'e}r{\^o}me}\binits{J.}} \AND
  \bauthor{\bsnm{Merlev{\`e}de},~\bfnm{Florence}\binits{F.}}
(\byear{2002}).
\btitle{Necessary and sufficient conditions for the conditional central limit
  theorem}.
\bjournal{Ann. Probab.}
\bvolume{30}
\bpages{1044--1081}.
\bid{doi={10.1214/aop/1029867121}, mr={1920101}}
\end{barticle}
\endbibitem

\bibitem[\protect\citeauthoryear{Fisher}{1929}]{fisher29}
\begin{barticle}[vtex]
\bauthor{\bsnm{Fisher},~\bfnm{R.~A.}\binits{R.~A.}}
(\byear{1929}).
\btitle{Tests of significance in harmonic analysis}.
\bjournal{Proc. Roy. Soc. Ser. A}
\bvolume{125}
\bpages{54--59}.
\end{barticle}
\endbibitem

%

\bibitem[\protect\citeauthoryear{Hall and Heyde}{1980}]{HH80}
\begin{bbook}[mr]
\bauthor{\bsnm{Hall},~\bfnm{P.}\binits{P.}} \AND
  \bauthor{\bsnm{Heyde},~\bfnm{C.~C.}\binits{C.~C.}}
(\byear{1980}).
\btitle{Martingale Limit Theory and Its Application: Probability and Mathematical Statistics}.
\bpublisher{Academic Press}, 
\baddress{New York}.
\bid{mr={624435}}
\end{bbook}
\endbibitem

%
%
%
\bibitem[\protect\citeauthoryear{Hunt and Young}{1974}]{HY74}
\begin{barticle}[mr]
\bauthor{\bsnm{Hunt},~\bfnm{Richard~A.}\binits{R.~A.}} \AND
  \bauthor{\bsnm{Young},~\bfnm{Wo~Sang}\binits{W.~S.}}
(\byear{1974}).
\btitle{A weighted norm inequality for {F}ourier series}.
\bjournal{Bull. Amer. Math. Soc.}
\bvolume{80}
\bpages{274--277}.
\bid{mr={0338655}}
\end{barticle}
\endbibitem

\bibitem[\protect\citeauthoryear{Kolmogorov}{1923}]{Kolmogorov23}
\begin{barticle}[vtex]
\bauthor{\bsnm{Kolmogorov},~\bfnm{A.}\binits{A.}}
(\byear{1923}).
\btitle{Une s\'erie de Fourier--Lebesgue divergente presque
partout}.
\bjournal{Fund. Math.}
\bvolume{4}
\bpages{324--328}.
\end{barticle}
\endbibitem

\bibitem[\protect\citeauthoryear{Lacey and Terwilleger}{2008}]{LT08}
\begin{barticle}[mr]
\bauthor{\bsnm{Lacey},~\bfnm{Michael}\binits{M.}} \AND
  \bauthor{\bsnm{Terwilleger},~\bfnm{Erin}\binits{E.}}
(\byear{2008}).
\btitle{A {W}iener--{W}intner theorem for the {H}ilbert transform}.
\bjournal{Ark. Mat.}
\bvolume{46}
\bpages{315--336}.
\bid{doi={10.1007/s11512-008-0080-2}, mr={2430729}}
\end{barticle}
\endbibitem

\bibitem[\protect\citeauthoryear{Lahiri}{2003}]{Lahiri03}
\begin{barticle}[mr]
\bauthor{\bsnm{Lahiri},~\bfnm{S.~N.}\binits{S.~N.}}
(\byear{2003}).
\btitle{A necessary and sufficient condition for asymptotic independence of
  discrete {F}ourier transforms under short- and long-range dependence}.
\bjournal{Ann. Statist.}
\bvolume{31}
\bpages{613--641}.
\bid{doi={10.1214/aos/1051027883}, mr={1983544}}
\end{barticle}
\endbibitem

\bibitem[\protect\citeauthoryear{Lin and Liu}{2009}]{LiuLin09}
\begin{barticle}[mr]
\bauthor{\bsnm{Lin},~\bfnm{Zhengyan}\binits{Z.}} \AND
  \bauthor{\bsnm{Liu},~\bfnm{Weidong}\binits{W.}}
(\byear{2009}).
\btitle{On maxima of periodograms of stationary processes}.
\bjournal{Ann. Statist.}
\bvolume{37}
\bpages{2676--2695}.
\bid{doi={10.1214/08-AOS590}, mr={2541443}}
\end{barticle}
\endbibitem

\bibitem[\protect\citeauthoryear{Olshen}{1967}]{Olshen67}
\begin{barticle}[mr]
\bauthor{\bsnm{Olshen},~\bfnm{Richard~A.}\binits{R.~A.}}
(\byear{1967}).
\btitle{Asymptotic properties of the periodogram of a discrete stationary
  process}.
\bjournal{J. Appl. Probab.}
\bvolume{4}
\bpages{508--528}.
\bid{mr={0228059}}
\end{barticle}
\endbibitem

\bibitem[\protect\citeauthoryear{Rootz{\'e}n}{1976}]{R76}
\begin{barticle}[mr]
\bauthor{\bsnm{Rootz{\'e}n},~\bfnm{Holger}\binits{H.}}
(\byear{1976}).
\btitle{Gordin's theorem and the periodogram}.
\bjournal{J. Appl. Probab.}
\bvolume{13}
\bpages{365--370}.
\bid{mr={0410876}}
\end{barticle}
\endbibitem

\bibitem[\protect\citeauthoryear{Rosenblatt}{1981}]{Rosenblatt81}
\begin{barticle}[mr]
\bauthor{\bsnm{Rosenblatt},~\bfnm{M.}\binits{M.}}
(\byear{1981}).
\btitle{Limit theorems for {F}ourier transforms of functionals of {G}aussian
  sequences}.
\bjournal{Z. Wahrsch. Verw. Gebiete}
\bvolume{55}
\bpages{123--132}.
\bid{doi={10.1007/BF00535155}, mr={608012}}
\end{barticle}
\endbibitem

\bibitem[\protect\citeauthoryear{Rosenblatt}{1985}]{Rosenblatt85}
\begin{bbook}[mr]
\bauthor{\bsnm{Rosenblatt},~\bfnm{Murray}\binits{M.}}
(\byear{1985}).
\btitle{Stationary Sequences and Random Fields}.
\bpublisher{Birkh\"auser}, \baddress{Boston, MA}.
\bid{mr={885090}}
\end{bbook}
\endbibitem

\bibitem[\protect\citeauthoryear{Schuster}{1898}]{Schuster1898}
\begin{barticle}[vtex]
\bauthor{\bsnm{Schuster},~\bfnm{A.}\binits{A.}}
(\byear{1898}).
\btitle{On the investigation of hidden periodicities with
  application to a supposed 26 day period of meteorological phenomena}.
\bjournal{Terrestrial Magnetism and Atmospheric Electricity}
\bvolume{3}
\bpages{13--41}.
\end{barticle}
\endbibitem

\bibitem[\protect\citeauthoryear{Terrin and Hurvich}{1994}]{TH94}
\begin{barticle}[mr]
\bauthor{\bsnm{Terrin},~\bfnm{Norma}\binits{N.}} \AND
  \bauthor{\bsnm{Hurvich},~\bfnm{Clifford~M.}\binits{C.~M.}}
(\byear{1994}).
\btitle{An asymptotic {W}iener--{I}t\^o representation for the low frequency
  ordinates of the periodogram of a long memory time series}.
\bjournal{Stochastic Process. Appl.}
\bvolume{54}
\bpages{297--307}.
\bid{doi={10.1016/0304-4149(94)00020-4}, mr={1307342}}
\end{barticle}
\endbibitem

\bibitem[\protect\citeauthoryear{Walker}{1965}]{walker65}
\begin{barticle}[mr]
\bauthor{\bsnm{Walker},~\bfnm{A.~M.}\binits{A.~M.}}
(\byear{1965}).
\btitle{Some asymptotic results for the periodogram of a stationary time
  series}.
\bjournal{J.~Aust. Math. Soc.}
\bvolume{5}
\bpages{107--128}.
\bid{mr={0177457}}
\end{barticle}
\endbibitem

\bibitem[\protect\citeauthoryear{Walker}{2000}]{walker00}
\begin{barticle}[mr]
\bauthor{\bsnm{Walker},~\bfnm{A.~M.}\binits{A.~M.}}
(\byear{2000}).
\btitle{Some results concerning the asymptotic distribution of sample {F}ourier
  transforms and periodograms for a discrete-time stationary process with a
  continuous spectrum}.
\bjournal{J. Time Ser. Anal.}
\bvolume{21}
\bpages{95--109}.
\bid{doi={10.1111/1467-9892.00175}, mr={1766176}}
\end{barticle}
\endbibitem

\bibitem[\protect\citeauthoryear{Woodroofe}{1992}]{Woodroofe92}
\begin{barticle}[mr]
\bauthor{\bsnm{Woodroofe},~\bfnm{Michael}\binits{M.}}
(\byear{1992}).
\btitle{A central limit theorem for functions of a {M}arkov chain with
  applications to shifts}.
\bjournal{Stochastic Process. Appl.}
\bvolume{41}
\bpages{33--44}.
\bid{doi={10.1016/0304-4149(92)90145-G}, mr={1162717}}
\end{barticle}
\endbibitem

\bibitem[\protect\citeauthoryear{Wu and Woodroofe}{2000}]{WW00}
\begin{barticle}[mr]
\bauthor{\bsnm{Wu},~\bfnm{Wei~Biao}\binits{W.~B.}} \AND
  \bauthor{\bsnm{Woodroofe},~\bfnm{Michael}\binits{M.}}
(\byear{2000}).
\btitle{A central limit theorem for iterated random functions}.
\bjournal{J.~Appl. Probab.}
\bvolume{37}
\bpages{748--755}.
\bid{mr={1782450}}
\end{barticle}
\endbibitem

\bibitem[\protect\citeauthoryear{Wu}{2005}]{wu05}
\begin{barticle}[mr]
\bauthor{\bsnm{Wu},~\bfnm{Wei~Biao}\binits{W.~B.}}
(\byear{2005}).
\btitle{Fourier transforms of stationary processes}.
\bjournal{Proc. Amer. Math. Soc.}
\bvolume{133}
\bpages{285--293}.
\bid{doi={10.1090/S0002-9939-04-07528-8}, mr={2086221}}
\end{barticle}
\endbibitem

\bibitem[\protect\citeauthoryear{Yajima}{1989}]{Yajima89}
\begin{barticle}[mr]
\bauthor{\bsnm{Yajima},~\bfnm{Yoshihiro}\binits{Y.}}
(\byear{1989}).
\btitle{A central limit theorem of {F}ourier transforms of strongly dependent
  stationary processes}.
\bjournal{J. Time Ser. Anal.}
\bvolume{10}
\bpages{375--383}.
\bid{doi={10.1111/j.1467-9892.1989.tb00036.x}, mr={1038470}}
\end{barticle}
\endbibitem

\bibitem[\protect\citeauthoryear{Wiener and Wintner}{1941}]{WW41}
\begin{barticle}[mr]
\bauthor{\bsnm{Wiener},~\bfnm{Norbert}\binits{N.}} \AND
  \bauthor{\bsnm{Wintner},~\bfnm{Aurel}\binits{A.}}
(\byear{1941}).
\btitle{On the ergodic dynamics of almost periodic systems}.
\bjournal{Amer. J. Math.}
\bvolume{63}
\bpages{794--824}.
\bid{mr={0006618}}
\end{barticle}
\endbibitem

\bibitem[\protect\citeauthoryear{Zygmund}{2002}]{Zygmund02}
\begin{bbook}[mr]
\bauthor{\bsnm{Zygmund},~\bfnm{A.}\binits{A.}}
(\byear{2002}).
\btitle{Trigonometric Series}. 
\bpublisher{Cambridge Univ. Press}, \baddress{Cambridge}.
\bid{mr={1963498}}
\end{bbook}
\endbibitem

\end{thebibliography}
\end{document}